\documentclass[11pt]{article}
\usepackage{amsmath,amssymb,theorem,amstext,amsgen,amsbsy,amsopn,amsfonts,graphicx,cases}

\usepackage{graphicx}
\usepackage{subfigure}
\usepackage{psfrag}
\usepackage{color}
\usepackage{enumerate}
\usepackage{epstopdf}
\usepackage[numbers,sort&compress]{natbib}
\usepackage{array}

\usepackage{latexsym}
\newtheorem{thm}{Theorem}[section]
\newtheorem{conj}[thm]{Conjecture}
\newtheorem{cor}[thm]{Corollary}
\newtheorem{lem}[thm]{Lemma}
\theorembodyfont{\rmfamily}

\def\pf{\bigskip\noindent {\bf Proof.}}

\def\dfn#1{{\sl #1}}

\def\less{\setminus}

\def\pf{\bigskip\noindent {\emph{Proof.}}}

\def\qed{ \hfill\vrule height3pt width6pt depth2pt}

\def\pf{\bigskip\noindent {\bf Proof.  }}

\title{An improved lower bound for the planar Tur\'an number of cycles}

\author{Yongxin Lan\thanks{School of Science, Hebei University of Technology, Tianjin 300401,   China.    Supported by the National Natural Science Foundation of China under Grant Nos. 12001154, 12161141006, Natural Science Foundation of Hebei Province under Grant No. A2021202025.    E-mail: yxlan@hebut.edu.cn}\hskip 1cm  Zi-Xia Song\thanks{Department  of Mathematics, University of Central Florida, Orlando, FL 32816, USA. Supported by the National   Science  Foundation under Grant No. DMS-1854903. E-mail:  Zixia.Song@ucf.edu }
}

\date{}

\begin{document}
\maketitle
\begin{abstract}
 The planar Tur\'an number of a graph $H$, denoted by $ex_{_\mathcal{P}}(n,H)$,   is the largest  number of edges in a  planar graph on $n $ vertices without containing $H$ as a subgraph. In this paper, we continue to study the topic of ``extremal" planar graphs initiated by Dowden [J. Graph Theory  83 (2016) 213--230].   We first obtain an improved lower bound for $ex_{_\mathcal{P}}(n,C_k)$ for all $k\ge 13$ and $n\ge 5(k-6+\lfloor{(k-1)}/2\rfloor)(k-1)/2$;  the construction for each  $k$ and $n$   provides a  simpler counterexample to a conjecture of Ghosh, Gy\H{o}ri, Martin, Paulos and Xiao [arxiv:2004.14094v1], which has recently been disproved   by  Cranston, Lidick\'y, Liu and Shantanam [Electron. J. Combin.  29(3) (2022) \#P3.31] for every $k\ge 11$ and $n$ sufficiently large (as a function of $k$).    We  then prove that $ex_{_\mathcal{P}}(n,H^+)=ex_{_\mathcal{P}}(n,H)$ for all $k\ge 5$ and $n\ge |H|+1$, where   $H\in\{C_k, 2C_k\}$ and $H^+$ is obtained from $H$ by adding a pendant edge to a vertex of degree two.

\end{abstract}

  {\bf Key words}.    planar Tur\'an number,  extremal planar graph, cycles, Theta graphs

  {\bf AMS subject classifications}.  05C10, 05C35

\baselineskip 16pt
\section{Introduction}

All graphs considered in this paper are finite and simple.  We use $K_n$,    $C_n$, $P_n$ and $K_{1, n-1}$ to denote the complete graph, cycle,   path, and star  on $n$ vertices, respectively.   For a positive integer $t$ and a graph $H$, we use $tH$ to denote the disjoint union of $t$ copies of   $H$, $\overline{H}$ the complement of $H$, and $H^+$ the graph obtained from $H$   by adding exactly one pendant edge to a vertex of degree two.
For an integer $k\ge 4$,  a \dfn{Theta graph} on $k$ vertices is the graph obtained from $C_k$ by adding a chord of the cycle.  Let $\Theta_k$ denote the family of all non-isomorphic Theta graphs on $k$ vertices and  let $\Theta_k^+=\{H^+\mid H\in \Theta_k\}$. Given a family $\mathcal{F}$, we say a graph is $\mathcal{F}$-free if it does not contain any graph in $\mathcal{F}$ as a subgraph.  When $\mathcal{F}=\{H\}$ we simply write $H$-free.  One of the best known results in extremal graph theory is Tur\'an's Theorem~\cite{Turan1941}, which
gives the maximum number of edges that a $K_k$-free graph on $n$ vertices can have. The celebrated  Erd\H{o}s-Stone Theorem~\cite{ErdosStone}  then extends this to the case when $K_k$ is
replaced by an arbitrary graph $H$, showing that the maximum number of edges possible
is $(1+o(1)){n\choose 2}\left(\frac{\chi(H)-2}{\chi(H)-1}\right)$, where $\chi(H)$ denotes the chromatic number of $H$.    Tur\'an-type problems when host graphs are hypergraphs are notoriously difficult.  Dowden~\cite{Dowden2016} in 2016 considered Tur\'an-type problems when host graphs are planar graphs, i.e.,  how many edges can an $\mathcal{F}$-free  planar graph on $n$ vertices have?  The \dfn{planar Tur\'an number of  $\mathcal{F}$}, denoted by $ex_{_\mathcal{P}}(n,\mathcal{F})$, is  the maximum number of edges in an $\mathcal{F}$-free  planar graph on $n$ vertices.   When $\mathcal{F}=\{H\}$ we write $ex_{_\mathcal{P}}(n,H)$.
 Dowden~\cite{Dowden2016}  obtained  tight bounds for  $ex_{_\mathcal{P}}(n, H) $ when $H\in\{C_4, C_5\}$; and shortly after,  this parameter has been investigated for various graphs $H$ in \cite{LZW, C6, Theta6, doublestar, 19LSSan,LSS,19LSS, 22LSS}. We refer the reader to a recent survey of Shi and the present   authors \cite{LSS21} for more information on this topic.
 \medskip

 It is worth noting that a  variety of sufficient conditions on   planar graphs $H$ such  that $ex_{_\mathcal{P}}(n,H)=3n-6$ for all $n\ge|H|$ has been investigated by Shi and  the present authors in \cite{19LSS}; in particular, $ex_{_\mathcal{P}}(n,H)=3n-6$ if $H$ contains three vertex-disjoint cycles or $H$ has at least seven vertices of degree at least four.     However,  determining  the  values of $ex_{_\mathcal{P}}(n,H)$,   when $H$ is a    subcubic graph without containing three vertex disjoint cycles,  remains wide open.
Very recently, Shi and the present authors~\cite{22LSS} completely determined the value of  $ex_{_\mathcal{P}}(n,H)$ when $H$ is cubic,   $H=C_3^+$, $H=2C_3$,  or  $H=C_3\cup C_3^+$.
In this paper, we continue to study $ex_{_\mathcal{P}}(n,H)$  and $ex_{_\mathcal{P}}(n,H^+)$ when $H\in \{C_k, 2C_k\}$; and     $ex_{_\mathcal{P}}(n,\Theta_k)$  and $ex_{_\mathcal{P}}(n,\Theta_k^+)$.  We first list some known results on $C_k$ and $\Theta_k$ that will be needed later.

 \begin{thm}[Dowden~\cite{Dowden2016}]\label{Dowden2016}  Let $n$ be a positive integer.
\begin{enumerate}[(a)]
\item $ex_{_\mathcal{P}}(n, C_3)=2n-4$ for all $n\geq 3$.
\item $ex_{_\mathcal{P}}(n, C_4)\leq {15}(n-2)/7$ for all $n\geq 4$, with equality when $n\equiv 30 (\rm{mod}\, 70)$.
\item  $ex_{_\mathcal{P}}(n, C_5)\leq (12n-33)/{5}$ for all $n\geq 11$. Equality holds for infinity many $n$.
\end{enumerate}
\end{thm}
\begin{thm}[Lan, Shi and Song~\cite{LSS}]\label{Theta4}
Let $n$ be a positive integer.
\begin{enumerate}[(a)]
\item   $ex_{_\mathcal{P}}(n,\Theta_4)\leq  {12(n-2)}/5$ for all $n\ge4$, with  equality when $n\equiv 12 (\rm{mod}\, 20)$.
\item $ex_{_\mathcal{P}}(n,\Theta_5)\leq {5(n-2)}/2$ for all $n\ge5$, with  equality when $n\equiv 50 (\rm{mod}\, 120)$.
\item $ex_{_\mathcal{P}}(n, C_6)\leq ex_{_\mathcal{P}}(n, \Theta_6) \le {18(n-2)}/7$ for all $n\geq6$.
\end{enumerate}
\end{thm}

Theorem~\ref{Theta4}(c) has been strengthened by the authors in \cite{C6, Theta6} with tight upper bounds.
 \begin{thm}[Ghosh et al.~\cite{C6, Theta6}]\label{C6Theta6}  Let $n$ be a positive integer.
\begin{enumerate}[(a)]

\item $ex_{_\mathcal{P}}(n, C_6)\leq (5n-14)/2$ for all $n\geq 18$, with  equality when $n\equiv 2 (\rm{mod}\, 5)$.

\item $ex_{_\mathcal{P}}(n, \Theta_6)\leq (18n-48)/7$ for all $n\geq 14$. Equality  holds for infinitely many $n$.

\end{enumerate}
\end{thm}

  In the same paper, Ghosh, Gy\H{o}ri, Martin, Paulo, and Xiao~\cite[Conjecture 15]{C6} conjectured that  \[ex_{_\mathcal{P}}(n, C_k)\le \left(3-\frac3k\right)n-6-\frac6k\] for all $k\ge 7$ and each sufficiently large $n$.  This has recently been disproved by Cranston, Lidick\'y, Liu and Shantanam~\cite{CLLS21} for every $k\ge 11$ and $n$ sufficiently large (as a function of $k$).

\begin{thm}[Cranston, Lidick\'y, Liu and Shantanam~\cite{CLLS21}]
Let $n$, $k$ and $\ell$ be positive integers with $k\ge11$. Let $\varepsilon$ be the remainder of $k$ when divided by $2$, and $\ell$ be even. Then for all
$n=((k-1){(5\ell-2)}/{2}+2)(\lfloor{3(k-1)}/{2}\rfloor-5)-(5\ell-4)$, we have $$ex_{_\mathcal{P}}(n,C_k)\ge \left(3-\frac{3}{(k-6+\lfloor\frac{k-1}{2}\rfloor)\frac{k-1}{k-3}-\frac{2}{k-3}}\right)n-6-\frac{12-\frac{42-6\varepsilon}{k-3}}{(k-6+\lfloor\frac{k-1}{2}\rfloor)\frac{k-1}{k-3}-\frac{2}{k-3}}.$$
\end{thm}

  They further proposed the following conjecture.

 \begin{conj}[Cranston, Lidick\'y, Liu and Shantanam~\cite{CLLS21}]
There exists a constant $D$ such that for all $k$ and for all sufficiently large $n$, we have
$$ex_{_{\mathcal{P}}}(n,C_k)\le \left(3-\frac{3}{Dk^{lg_2^3}}\right)n.$$
\end{conj}

 In Section~\ref{s:lbd}, we establish an improved lower bound for $ex_{_{\mathcal{P}}}(n,C_k)$ by proving that $ex_{_\mathcal{P}}(n, C_k)> \left(3-3/k\right)n-6-6/k$ for all $k\ge 13$ and   $n\ge 5(k-6+\lfloor{(k-1)}/2\rfloor)(k-1)/2$;   our  construction then provides a  simpler counterexample to the conjecture of Ghosh et al. \cite{C6} for all $k\ge 13$ and $n\ge 5(k-6+\lfloor{(k-1)}/2\rfloor)(k-1)/2$.  Using a recursive relation, we then prove that $ex_{_{\mathcal{P}}}(n,C_k)=ex_{_{\mathcal{P}}}(n,C_k^+)$ for all $k\ge 4$ and $n\ge k+1$, except when $n=8$ and $k=4$ in Section~\ref{s:Ck+};  $ex_{_{\mathcal{P}}}(n,2C_k)=ex_{_{\mathcal{P}}}(n,C_k^+\cup C_k)$ for all $k\ge 4$ and $n\ge 2k+1$  in Section~\ref{s:2Ck}; and   $ex_{_{\mathcal{P}}}(n,\Theta_k)=ex_{_{\mathcal{P}}}(n,\Theta_k^+)$ for all $k\ge 4$ and $n\ge k+2$ in Section~\ref{s:Theta+}.   \medskip

We end this section by introducing more notation.
 Given a graph $G$, we will use $V(G)$ to denote the vertex set, $E(G)$ the edge set, $|G|$ the number of vertices, $e(G)$ the number of edges, $\delta(G)$ the minimum degree, $\Delta(G)$ the maximum degree. For a vertex $v\in V(G)$, we will use $N_G(v)$ to denote the set of vertices in $G$ which are adjacent to $v$. Let $d_G(v)=|N_G(v)|$ denote the degree of the vertex $v$ in $G$   and $N_{G}[v]=N_{G}(v)\cup \{v\}$.
 For any set  $S \subset V(G)$,
the subgraph of $G$ induced on  $S$, denoted by $G[S]$, is the graph with vertex set $S$ and edge set $\{xy \in E(G) \mid  x, y \in S\}$. We denote by   $G \less S$ the subgraph of $G$ induced  on
$V(G) \less S$.   If $S=\{v\}$, then we simply write   $G\less v$.
The {\dfn{join}} $G+H$ (resp. {\dfn{union}} $G\cup H$)
of two vertex-disjoint graphs
$G$ and $H$ is the graph having vertex set $V(G)\cup V(H)$  and edge set $E(G)
\cup E(H)\cup \{xy\, |\,  x\in V(G),  y\in V(H)\}$  (resp. $E(G)\cup E(H)$). 
Given two isomorphic graphs $G$ and $H$,   we may (with a slight but common abuse of notation) write $G = H$.  For any positive integer $k$,  we define   $[k]:=\{1,2, \ldots, k\}$.

\section{An improved lower bound for $ex_{_{\mathcal{P}}}(n,C_k)$}\label{s:lbd}
\begin{figure}[htbp]
  \centering
  \includegraphics[scale=0.25]{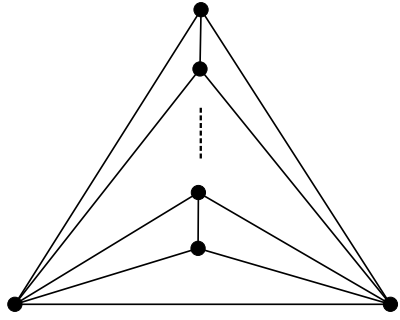}
  \caption{The construction of $\mathcal{T}_m$.}\label{Tm}
\end{figure}

\begin{figure}[htbp]
  \centering
  \includegraphics[scale=0.14]{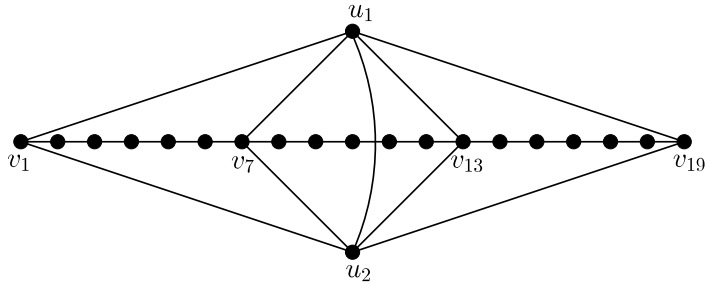}
  \caption{The construction of $R_\ell$ when $\ell=21$ and $k=7$.}\label{R}
\end{figure}

 In this section we construct a $C_k$-free planar graph $G$ on $n\ge k$ vertices  with $e(G)=3n-6$ for all  $k\ge 11$ and  $n\le (3k-11)/2$, and  $e(G)\ge\left(3-{(3-\frac{2}{k-1})}/{(k-6+\lfloor{(k-1)}/2\rfloor)}\right)n-11> \left(3-3/k\right)n-6-6/k$   for all $k\ge 13$ and $n\ge 5(k-6+\lfloor{(k-1)}/2\rfloor)(k-1)/2$.  It follows that our construction $G$ provides a  simpler counterexample to the conjecture of Ghosh et al. \cite[Conjecture 15]{C6} for all $k\ge 13$ and $n\ge 5(k-6+\lfloor{(k-1)}/2\rfloor)(k-1)/2$. 
\begin{thm}\label{LB-Ck}
Let $n$ and $k$ be  integers with $n\ge k\ge 11$. Let $r$   be the remainder of  $n-4$ when divided by $k-6+\lfloor(k-1)/2\rfloor$. Then
\begin{enumerate}[($i$)]
  \item $ex_{_{\mathcal{P}}}(n,C_k)=3n-6$ for all $n\le  k-5+\left\lfloor\frac{k-1}{2}\right\rfloor$, and
  \item $ex_{_{\mathcal{P}}}(n,C_k)\ge\left(3-\frac{3-\frac{2}{k-1}}{k-6+\left\lfloor\frac{k-1}2\right\rfloor}\right)n+
  \frac{12+3r-\frac{8+2r}{k-1}}{k-6+\left\lfloor\frac{k-1}2\right\rfloor}+\frac{4}{k-1}-\min\{r+10,11\}$
   for all  $n\ge  k-4+\left\lfloor\frac{k-1}{2}\right\rfloor$.
\end{enumerate}
\end{thm}
\pf Let $n$, $k$ and $r $  be given as in the statement. Throughout the proof,  let $\mathcal{T}_m:=K_2+P_{m-2}$ be the plane triangulation on $m$ vertices,  where $m\ge3$ is an integer and   $\mathcal{T}_m$ is depicted in  Figure~\ref{Tm}. 
Note that  $\mathcal{T}_m$ has exactly $2m-4$ $3$-faces.
For any integer $\ell\ge4$, let $a:=\lfloor {(\ell-3)}/{(k-1)}\rfloor$. 
Let $R_\ell$ be a planar graph on $\ell$ vertices obtained from a path $P:=v_1v_2\ldots v_{\ell-2}$ by adding two new adjacent vertices, say $u_1,u_2$, such that for each $i\in[2]$, $u_i$ is adjacent to $v_1,v_{k-1+1},v_{2(k-1)+1},\ldots,v_{a(k-1)+1}$,
which is depicted in Figure~\ref{R} when $\ell=21$ and $k=7$.  Note that $e(R_\ell)=(\ell-3)+2(a+1)+1=\ell+2\lfloor {(\ell-3)}/{(k-1)}\rfloor$. We proceed the proof by considering the parity of $k$.

We first consider the case when $k\ge 11$ is odd. Let $k:=2p+1$, where $p\ge5$ is an integer. For each $n$ satisfying $2p+1=k\le n\le  k-5+\lfloor(k-1)/2\rfloor=3p-4=p+(2p-4)$, let $\mathcal{T}_n^p$ denote a plane  triangulation on $n$ vertices obtained from $\mathcal{T}_p$ by adding $n-p\le 2p-4$ new vertices: each to a $3$-face $F$ of $\mathcal{T}_p$ and then joining it to all vertices on the boundary of $F$. Note that $\mathcal{T}_n^{p}\less V(\mathcal{T}_p)$ contains no edge and each cycle in $\mathcal{T}_n^{p}$ contains at most $p$ vertices of $\mathcal{T}_n^{p}\less V(\mathcal{T}_p)$. It follows that $\mathcal{T}_n^{p}$ is $C_k$-free, and so $ex_{_{\mathcal{P}}}(n,C_k)=e(\mathcal{T}_n^{p})=3n-6$ for all $n$ satisfying $k\le n\le  k-5+\lfloor(k-1)/2\rfloor$, as desired.
We may assume that $n\ge   k-4+\lfloor(k-1)/2\rfloor=3p-3$.
Let $t\ge0$ be an integer satisfying  $$t(3p-5)+r=n-4.$$
By the choice of $r$, we have $0\le r\le k-7+\lfloor(k-1)/2\rfloor=3p-6$.

Let $H_i:=\mathcal{T}_{3p-4}^{p}$ for each $i\in[t]$;
 $H_{t+1}$ be the complete graph $K_{r+2}$ when $r\le 1$, the plane triangulation $\mathcal{T}_{r+2}$ when $2\le r\le k-3$, the plane triangulation $\mathcal{T}_{r+2}^{p}$ when $k-2\le r\le 3p-6$. For each $j\in[t+1]$, let $y$ and $z$ be any two distinct vertices on the outer-face of $H_j$; it is easy to check that each cycle of $H_j$ has  length at most $k-1$.
Let $G$ be the planar graph obtained from $R_{t+4}$ by identifying each edge $v_iv_{i+1}$ with the edge $yz$ in $H_i$ for any $i\in[t+1]$.
Then $G$ is $C_k$-free because any cycle containing $u_1$ or $u_2$ has length either at least $k+1$ or at most $4$.  Note that $|G|=t+4+t(3p-6)+r$, and  $|H_{t+1}|=2$ with equality when $r=0$.
Therefore,
\begin{align*}
ex_{_\mathcal{P}}(n,C_k)\ge e(G)
&=e(R_{t+4})+(e(H_1)-1)+(e(H_2)-1)+\cdots+(e(H_{t+1})-1)\\
&=t+4+2\left\lfloor \frac{t+1}{k-1}\right\rfloor+t(3(3p-4)-7)+ \max\{3(r+2)-7,0\}\\
&=3n-3t+2\left\lfloor \frac{t+1}{k-1}\right\rfloor-\min\{r+8,9\}\\
&\ge 3n-3t+2\cdot\frac{t+1-(k-2)}{k-1}-\min\{r+8,9\}\\
&=3n-\left(3-\frac{2}{k-1}\right)t-\frac{2k-6}{k-1}-\min\{r+8,9\}\\
&=3n-\frac{\left(3-\frac{2}{k-1}\right)(n-4-r)}{3p-5}+\frac{4}{k-1}-\min\{r+10,11\}\\
&=\left(3-\frac{3-\frac{2}{k-1}}{k-6+\frac{k-1}2}\right)n+\frac{\left(3-\frac{2}{k-1}\right)(4+r)}{k-6+\frac{k-1}2}+\frac{4}{k-1}-\min\{r+10,11\}\\
&= \left(3-\frac{3-\frac{2}{k-1}}{k-6+\frac{k-1}2}\right)n+
  \frac{12+3r-\frac{8+2r}{k-1}}{k-6+\frac{k-1}2}+\frac{4}{k-1}-\min\{r+10,11\}.
 \end{align*}

We next consider the case when $k\ge12$ is even.  Let $k:=2p$, where $p\ge6$ is an integer. For each $n$ satisfying $2p=k\le n\le k-5+\lfloor(k-1)/2\rfloor=3p-6=(p-1)+(2p-5)$, let $\mathcal{L}_n^{p-1}$ denote the plane  triangulation on $n$ vertices obtained from $\mathcal{T}_{p-1}$ by first adding   two adjacent new vertices to a $3$-face $F'$ of $\mathcal{T}_{p-1}$ with one   vertex joining to all vertices on the boundary of $F'$, the other to exactly two vertices on the boundary of $F'$; then adding $n-p-1 \le 2p-7$ new vertices: each to a $3$-face $F\ne F'$ of $\mathcal{T}_{p-1}$ and then joining it to all vertices on the boundary of $F$. Note that  $\mathcal{L}_n^{p-1}\less V(\mathcal{T}_{p-1})$ contains exactly one edge and each cycle in $\mathcal{L}_n^{p-1}$ contains at most $p$ vertices of  $\mathcal{L}_n^{p-1}\less V(\mathcal{T}_{p-1})$. Hence $\mathcal{L}_n^{p-1}$ is $C_k$-free. It follows that
$ex_{_{\mathcal{P}}}(n,C_k)=e(\mathcal{L}_n^{p-1})=3n-6$ for all $n\le k-5+\lfloor(k-1)/2\rfloor$, as desired.
We may assume that $n\ge   k-4+\lfloor(k-1)/2\rfloor=3p-5$.
Let $t\ge0$ be an integer satisfying  $$t(3p-7)+r=n-4.$$
By the choice of $r$, we have $0\le r\le k-7+\lfloor(k-1)/2\rfloor=3p-8$.

Let $H_i:=\mathcal{L}_{3p-6}^{p-1}$ for each $i\in[t]$;
 $H_{t+1}$ be the complete graph $K_{r+2}$  when $r\le1$, the plane triangulation $\mathcal{T}_{r+2}$ when $2\le r\le k-3$, the plane triangulation $\mathcal{L}_{r+2}^{p-1}$ when $k-2\le r\le 3p-8$. For each $j\in[t+1]$, let $y$ and $z$ be any two distinct vertices on the outer-face of $H_j$; it is easy to check that each cycle of $H_j$ has  length at most $k-1$.
Let $G$ be the planar graph obtained from $R_{t+4}$ by identifying each edge $v_iv_{i+1}$ with the edge $yz$ in $H_i$ for any $i\in[t+1]$.
It follows that $G$ is $C_k$-free and $|G|=t+4+t(3p-8)+r$. Note that $|H_{t+1}|=2$ with equality when $r=0$.
Therefore,
\begin{align*}
ex_{_\mathcal{P}}(n,C_k)\ge e(G)
&=e(R_{t+4})+(e(H_1)-1)+(e(H_2)-1)+\cdots+(e(H_{t+1})-1)\\
&=t+4+2\left\lfloor \frac{t+1}{k-1}\right\rfloor+t(3(3p-6)-7)+ \max\{3(r+2)-7,0\}\\
&=3n-3t+2\left\lfloor \frac{t+1}{k-1}\right\rfloor-\min\{r+8,9\}\\
&\ge 3n-3t+2\cdot\frac{t+1-(k-2)}{k-1}-\min\{r+8,9\}\\
&=3n-\left(3-\frac{2}{k-1}\right)t-\frac{2k-6}{k-1}-\min\{r+8,9\}\\
&=3n-\frac{\left(3-\frac{2}{k-1}\right)(n-4-r)}{3p-7}+\frac{4}{k-1}-\min\{r+10,11\}\\
&=\left(3-\frac{3-\frac{2}{k-1}}{k-6+\left\lfloor\frac{k-1}2\right\rfloor}\right)n+\frac{\left(3-\frac{2}{k-1}\right)(4+r)}{k-6+\left\lfloor\frac{k-1}2\right\rfloor}+\frac{4}{k-1}-\min\{r+10,11\}.\\
&= \left(3-\frac{3-\frac{2}{k-1}}{k-6+\left\lfloor\frac{k-1}2\right\rfloor}\right)n+
  \frac{12+3r-\frac{8+2r}{k-1}}{k-6+\left\lfloor\frac{k-1}2\right\rfloor}+\frac{4}{k-1}-\min\{r+10,11\}.
 \end{align*}
This completes the proof of Theorem~\ref{LB-Ck}.\qed

\section{Planar Tur\'an numbers for  $C_k^+$ and $C_k$}\label{s:Ck+}
For the remainder of this paper,  let $O_p$ denote a maximal  outer-planar graph on $p\ge1$ 
vertices with maximum degree at most four. Note that $O_p$ has exactly $\max\{2p-3,0\}$ edges. 
It has recently been shown   in \cite{22LSS} that $ex_{_\mathcal{P}}(n,C_3^+)=   ex_{_\mathcal{P}}(n,C_3)$ for all $n\ge 4$. In this section we extend this result to all $k\ge 4$ by proving  that $ex_{_\mathcal{P}}(n,C_k^+)=   ex_{_\mathcal{P}}(n,C_k)$ for all $n\ge k+1$, except when $n=8$ and $k=4$.
We first consider the case when $k\ge 5$. We begin with a useful lemma which establishes  a  recursive relation for $ex_{_{\mathcal{P}}}(n,C_k)$ when $n\ge 2k+1\ge7$. The same idea will be explored in the next three sections.

\begin{lem}\label{Ck}
Let $n$ and $k$ be integers with $n\ge2k+1$ and $k \ge3$. Then \[ex_{_{\mathcal{P}}}(n,C_k)\ge ex_{_{\mathcal{P}}}(n-k,C_k)+3k-6.\]
\end{lem}

\pf   For integers $n$ and $k$   with $n\ge2k+1 \ge7$, it suffices to construct a $C_k$-free planar graph on $n $ vertices with $ex_{_{\mathcal{P}}}(n-k,C_k)+3k-6$ edges. Let $H$ be any $C_k$-free planar graph with  $|H|=n-k\ge k+1$ and  $e(H)=ex_{_{\mathcal{P}}}(n-k,C_k)$. Let $H^*$ be the disjoint union of $O_{k-2}$ and $K_2$ when $k\ge 4$,  and $H^*=\overline{K}_k$ when $k=3$. Finally, let $v\in V(H)$ and let $G$ be the planar graph obtained from the disjoint union of $H$ and $H^*$ by   joining $v$ to every  vertex in $H^*$.
It is simple to check that  $G$ is   $C_k$-free   with $|G|=n$ and
\begin{align*}
e(G) &= \begin{cases} e(H)+e(O_{k-2})+e(K_2)+k   &\text{\hskip 2.85cm  if } k\ge4\\
e(H)+k  &\text{\hskip 2.85cm  if } k=3\\
\end{cases}\\
&= \begin{cases} ex_{_{\mathcal{P}}}(n-k,C_k)+(2(k-2)-3)+1+k &\text{\hskip 1cm if } k\ge4\\
ex_{_{\mathcal{P}}}(n-k,C_k)+k &\text{\hskip 1cm  if } k=3\\
\end{cases}\\
&=ex_{_{\mathcal{P}}}(n-k,C_k)+3k-6.
\end{align*}
Hence,  $ex_{_{\mathcal{P}}}(n,C_k)\ge e(G)= ex_{_{\mathcal{P}}}(n-k,C_k)+3k-6$, as desired.\qed

\begin{thm}\label{Ck+}
Let $n$ and $k$ be integers with $n\ge k+1$ and $k\ge5$. Then \[ex_{_\mathcal{P}}(n,C_k^+)=ex_{_\mathcal{P}}(n,C_k).\]
\end{thm}
\pf Let $n$ and $k$ be integers with $n\ge k+1$ and $k\ge5$.  Since  $C_k$ is a subgraph of $C_k^+$, we see that   $ex_{_{\mathcal{P}}}(n,C_k^+)\ge ex_{_{\mathcal{P}}}(n,C_k)$. We next proceed the proof by induction on $n$ to   show  that $ ex_{_{\mathcal{P}}}(n,C_k^+)\le ex_{_{\mathcal{P}}}(n,C_k)$.  Let $G$ be a $C_k^+$-free planar graph on $n\ge k+1$ vertices with $ex_{_{\mathcal{P}}}(n,C_k^+)$ edges. If $G$ is $C_k$-free, then $e(G)\le ex_{_{\mathcal{P}}}(n,C_k)$, as desired. We may assume that $G$ contains $C_k$ as a subgraph. Let $A$ be the set of vertices of a copy of $C_k$ in $G$. Then $G[A]$ is a component of $G$ with $e(G[A])\le 3k-6$ because $G$ is $C_k^+$-free and $k\ge 5$. In addition, $G\less A$ is   a $C_k^+$-free planar graph on $n-k\ge 1$ vertices. 
\medskip

We first show that the statement holds  when $k+1\le n\le 2k$. Note that  $e(G\less A)\le {n-k\choose2}$ when $k+1\le n\le k+2$, and $e(G\less A)\le 3(n-k)-6$ when $k+3\le n\le 2k$. Hence,
\begin{eqnarray*}
  e(G)=e(G[A])+e(G\less A)&\le &
  \begin{cases}
  3k-6 & \text{if}\,\quad  n=k+1 \\
  3k-6+1 & \text{if}\,\quad  n=k+2 \\
  3k-6+3(n-k)-6 & \text{if}\,\quad k+3\le n\le 2k \\
\end{cases} \\
  &=&\begin{cases}
  3n-9 & \text{if}\,\quad  n=k+1 \\
  3n-11 & \text{if}\,\quad  n=k+2 \\
  3n-12 & \text{if}\,\quad k+3\le n\le 2k. \\
\end{cases}
\end{eqnarray*}
We then  construct a $C_k$-free planar graph $G^*$ on $n$ vertices with $e(G^*)\ge e(G)$.  Let $n=t(k-2)+\varepsilon+1$, where $t=1$ and $2\le\varepsilon=n-k+1\le k-2$ when $k+1\le n\le 2k-3$;   $t=2$ and $1\le \varepsilon=n-2k+3\le 3$ when $2k-2\le n\le 2k$. Let $G^*:=K_1+(tO_{k-2}\cup O_{\varepsilon})$.
Then   $G^*$ is a $C_k$-free planar graph with $|G^*| =n$  and
\begin{align*}
e(G^*)&= (t(k-2)+\varepsilon)+t(2(k-2)-3)+\max\{2\varepsilon-3,0\}\\
&=3(t(k-2)+\varepsilon)-2\varepsilon-3t+\max\{2\varepsilon-3,0\}\\
&=3(n-1)-3t+\max\{-3,-2\varepsilon\}\\
&=
  \begin{cases}
  3n-9 & \text{if}\,\quad  k+1\le n\le2k-3 \\
  3n-11 & \text{if}\,\quad  n=2k-2 \\
  3n-12 & \text{if}\,\quad 2k-1\le n\le 2k \\
\end{cases}\\
&\ge e(G).
\end{align*}
  Since $G^*$ is $C_k$-free, we see  that   $e(G)\le e(G^*)\le ex_{_{\mathcal{P}}}(n,C_k)$ for all $n$ satisfying  $k+1\le n\le 2k$.  \medskip

We may assume that $n\ge 2k+1$. Then $n-k\ge k+1$. By the induction hypothesis, we have $ex_{_{\mathcal{P}}}(n-k,C_k^+)\le ex_{_{\mathcal{P}}}(n-k,C_k)$. Note that  $e(G\less A)\le ex_{_{\mathcal{P}}}(n-k,C_k^+)$ because $G\less A$ is $C_k^+$-free. Hence,
\[e(G)=e(G[A])+e(G\less A)\le 3k-6+ex_{_{\mathcal{P}}}(n-k,C_k^+)\le3k-6+ex_{_{\mathcal{P}}}(n-k,C_k).\]
By Lemma~\ref{Ck}, we see that $e(G)\le ex_{_{\mathcal{P}}}(n,C_k)$, as desired.\medskip

This completes the proof of Theorem~\ref{Ck+}.\qed\bigskip

We end this section by considering the planar Tur\'an number for $C_4^+$. We need to determine the exact value of $ex_{_\mathcal{P}}(8,C_4)$.

\begin{lem}\label{n=8} $ex_{_\mathcal{P}}(8,C_4)=11$.
\end{lem}

\pf   Let $G:=C_8$  with vertices $v_1,\ldots,v_8$ in order. Let $G^*$ be the graph obtained from $G$ by adding edges $v_iv_{i+2}$ for each $i\in\{1,3,5\}$. Then  $G^*$ is a $C_4$-free planar graph with $|G^*|=8$ and $e(G^*)=11$. Hence, $ex_{_\mathcal{P}}(8,C_4)\ge e(G^*)=11$.   Next suppose   $ex_{_\mathcal{P}}(8,C_4)\ge 12$. Let $H$ be a $C_4$-free planar graph on $8$ vertices with $e(H)\ge12$. Let $v\in V(H)$ be  such that $d_H(v)=\Delta(H)$ and   let $B:=V(H)\less N_H[v]$. Then $d_H(v)\ge 3$, $|B|=7-d_H(v)\le 4$, and $e(H)= e(H[N_H[v]])+e_H(N_H(v),B)+e(H[B])$. In addition,  $H[N_H(v)]$ is $P_3$-free and each vertex in $B$ is adjacent to at most one vertex in $N_H(v)$. It follows that  $e(H[N_H(v)])\le \lfloor d_H(v)/2\rfloor$ and $e_H(N_H(v),B)\le |B|$. If $d_H(v)\ge4$, then $|B|\le 3$ and  $e(H[B])\le {|B|\choose 2}$.  Thus $12\le e(H)\le d_H(v)+\lfloor d_H(v)/2\rfloor+|B|+{|B|\choose 2}=7+\lfloor d_H(v)/2\rfloor+{7-d_H(v)\choose 2}\le 12$. It follows that  $d_H(v)=4$,  $e(H[N_H(v)])=\lfloor d_H(v)/2\rfloor=2$, $e_H(N_H(v),B)=|B|=3$ and $ H[B]=K_3$. One can check that $H$ contains  $C_4$ as a subgraph, a contradiction. Thus $d_H(v)=3$. Then $|B|=4$ and $e(H[B])\le 4$ because $H$ is $C_4$-free.  Thus $12\le e(H)\le d_H(v)+\lfloor d_H(v)/2\rfloor+|B|+4=12$. It follows that
$e_H(N_H(v),B)=|B|=4$ and $H[B]=K_1+({K_2\cup K_1})$. But then $\Delta(H)\ge4$, a contradiction.
 Thus, $ex_{_\mathcal{P}}(8,C_4)\le 11$ and so $ex_{_\mathcal{P}}(8,C_4)=11$.
\qed

\begin{thm}\label{C4+}
  Let $n\ge5$ be integer. Then
  \[ex_{_\mathcal{P}}(n,C_4^+)=
  \begin{cases}
    ex_{_\mathcal{P}}(n,C_4) & \mbox{if } n\ne8 \\
    ex_{_\mathcal{P}}(n,C_4)+1=12 & \mbox{if } n=8.
  \end{cases}  \]
\end{thm}
\pf  Since $C_4$ is a subgraph of $C_4^+$ and $2K_4$ is $C_4^+$-free, we see that  $ex_{_{\mathcal{P}}}(n,C_4^+)\ge ex_{_{\mathcal{P}}}(n,C_4)$,  and  $ex_{_\mathcal{P}}(8,C_4^+)\ge e(2K_4)=12=ex_{_\mathcal{P}}(8,C_4)+1$ by Lemma~\ref{n=8}.  We next   prove that $ ex_{_{\mathcal{P}}}(n,C_4^+)\le ex_{_{\mathcal{P}}}(n,C_4)$ when $n\ne8$ and $ex_{_{\mathcal{P}}}(n,C_4^+)\le 12=ex_{_{\mathcal{P}}}(n,C_4)+1$ when $n=8$. We shall proceed the proof by induction on $n$. Let $G$ be a $C_4^+$-free planar graph on $n$ vertices with $ex_{_{\mathcal{P}}}(n,C_4^+)$ edges. If $G$ is $C_4$-free, then $e(G)\le ex_{_{\mathcal{P}}}(n,C_4)$, as desired. We may assume that $G$ contains $C_4$ as a subgraph. Let $A$ be the set of vertices of a copy of $C_4$ in $G$. Then $G[A]$ is a component of $G$ because $G$ is $C_4^+$-free. Hence, $e(G)=e(G[A])+e(G\less A)\le 6+e(G\less A)$.\medskip

We first prove that the statement holds for  $5\le n\le 8$.  In this case, we have $e(G\less A)\le {n-4\choose2}$.  When $n=8$, we have $e(G)\le 6+{n-4\choose2}=12=ex_{_\mathcal{P}}(n,C_4)+1$. When $5\le n\le 7$,  note that $M_{n}: =K_1+(\lfloor (n-1)/2\rfloor K_2\cup (\lceil (n-1)/2\rceil-\lfloor (n-1)/2\rfloor)K_1)$ is $C_4$-free.  It follows that   $e(G)\le 6+{n-4\choose2}=e(M_n)\le ex_{_{\mathcal{P}}}(n,C_4)$. Hence the statement holds for  $5\le n\le 8$.\medskip

 We next prove that the statement holds for  $9\le n\le12$.    In this case, we have  $5\le |G\less A|\le 8$.  When $9\le n\le 11$, we have  $e(G\less A)\le  ex_{_{\mathcal{P}}}(n-4,C_4)$ and so $e(G)\le 6+ex_{_{\mathcal{P}}}(n-4,C_4)$. By Lemma~\ref{Ck}, we have $e(G)\le ex_{_{\mathcal{P}}}(n,C_4)$ for any $9\le n\le11$. When $n=12$, we see that $|G\less A|=8$ and so $e(G\less A)\le   ex_{_{\mathcal{P}}}(8,C_4^+)= ex_{_{\mathcal{P}}}(8,C_4)+1=12$. Hence     $e(G)\le 6+12=18$. We next construct a $C_4$-free planar graph with  $12$ vertices and $18$ edges.  Let $H:=C_{12}$   with vertices $v_1,\ldots,v_{12}$ in order. Let $H^*$ be the planar  graph obtained from $H$ by adding edge $v_iv_{i+2}$ for each $i\in\{1,3,5,7,9,11\}$, where the arithmetic on the index $i+2$ is done modulo $12$.  One can see that  $H^*$ is $C_4$-free   with $|H^*|=12$ and $e(H^*)=18$. Hence, $e(G)\le 18=e(H^*)\le ex_{_{\mathcal{P}}}(12,C_4)$. This proves that the statement holds for $9\le n\le12$.   \medskip

We may  assume that $n\ge13$. Then $|G\less A|\ge 9$ and $G\less A$ is $C_4^+$-free.  By the induction hypothesis, $e(G\less A)\le ex_{_{\mathcal{P}}}(n-4,C_4^+)\le ex_{_{\mathcal{P}}}(n-4,C_4)$. Hence,
$$e(G)\le 6+e(G\less A)\le 6+ex_{_{\mathcal{P}}}(n-4,C_4^+)\le6+ex_{_{\mathcal{P}}}(n-4,C_4).$$
By Lemma~\ref{Ck}, we see that $e(G)\le ex_{_{\mathcal{P}}}(n,C_4)$. This completes the proof of Theorem~\ref{C4+}.\qed\bigskip

Corollary~\ref{c:Ck+} below  follows directly from Theorem~\ref{Ck+}, Theorem~\ref{C4+}, Theorem~\ref{Dowden2016}($b,c$),  Theorem~\ref{C6Theta6}($a$) and Theorem~\ref{LB-Ck}.

\begin{cor}\label{c:Ck+}
Let $n, k$ be   positive integers. Then

 \begin{enumerate}[(a)]

\item $ex_{_\mathcal{P}}(n, C_4^+)\leq {15}(n-2)/7$ for all $n\geq 4$, with equality when $n\equiv 30 (\rm{mod}\, 70)$.

\item  $ex_{_\mathcal{P}}(n, C_5^+)\leq (12n-33)/{5}$ for all $n\geq 11$. Equality holds for infinity many $n$.

\item  $ex_{_\mathcal{P}}(n, C_6^+)\leq (5n-14)/2$ for all $n\geq 18$, with  equality when $n\equiv 2 (\rm{mod}\, 5)$.

\item $ex_{_{\mathcal{P}}}(n,C_k^+)=3n-6$ for all $k\ge 11$ and $k+1\le n\le  k-5+\lfloor(k-1)/2\rfloor$, and   \[ex_{_{\mathcal{P}}}(n,C_k^+)\ge \left(3-\frac{3-\frac{2}{k-1}}{k-6+\left\lfloor\frac{k-1}2\right\rfloor}\right)n+
  \frac{12+3r-\frac{8+2r}{k-1}}{k-6+\left\lfloor\frac{k-1}2\right\rfloor}+\frac{4}{k-1}-\min\{r+10,11\}\] for all $k\ge11$ and  $n\ge  k-4+\lfloor{(k-1)}/{2}\rfloor$, where $r$ is the remainder of  $n-4$ when divided by $k-6+\lfloor(k-1)/2\rfloor$.
\end{enumerate}
\end{cor}

\section{Planar Tur\'an numbers for  $C_k^+\cup C_k$ and $2C_k$}\label{s:2Ck}
It has recently been shown  in \cite{22LSS} that $ex_{_\mathcal{P}}(n,C_k\cup C_k^+)=ex_{_\mathcal{P}}(n,2C_k)$ for each  $k\in\{3,4,5,6\}$ and all $n\ge 2k+1$. Using a different method,   we extend this result further by proving that  $ex_{_\mathcal{P}}(n,C_k\cup C_k^+)=ex_{_\mathcal{P}}(n,2C_k)$ for all $n\ge 2k+1$ and $k\ge4$. We need a lemma from  \cite{22LSS}.

\begin{lem}[\cite{22LSS}]\label{22LSS}
Let $n$ and $k$ be positive integers with  $n\ge 2k\ge8$. Then $$ex_{_\mathcal{P}}(n,2C_k)\ge\left(3-\frac{1}{k-2}\right)n+\frac{3+r}{k-2} -5+\max\{1-r, 0\},$$ where $r$ is the remainder of $n-3$ when divided by $k-2$.
 \end{lem}

We first prove   recursive relations for $ex_{_{\mathcal{P}}}(n,2C_k)$ when  $n\ge 4k+1\ge17$.
\begin{lem}\label{2Ck}
Let $n$ and $k$ be integers with $n\ge4k+1$ and $k\ge4$. Then
\begin{enumerate}[(a)]
  \item $ex_{_{\mathcal{P}}}(n,2C_k)\ge ex_{_{\mathcal{P}}}(n-2k,2C_k)+6k-12$.
  \item $ex_{_{\mathcal{P}}}(n,2C_k)\ge ex_{_{\mathcal{P}}}(n-2k,C_k)+6k-6$.
\end{enumerate}

\end{lem}

\pf  To prove ($a$), it suffices to   construct a $2C_k$-free planar graph on $n$ vertices with $ex_{_{\mathcal{P}}}(n-2k,2C_k)+6k-12$ edges. Let $H_1$ be any $2C_k$-free planar graph  with  $|H_1|=n-2k\ge 2k+1$  and $e(H_1)=ex_{_{\mathcal{P}}}(n-2k,2C_k)$. Let $v_1\in V(H_1)$. Let $G_1$ be the planar graph obtained from the disjoint union of $H_1\cup 2O_{k-2}\cup 2K_2$ by joining $v_1$ to every vertex of $2O_{k-2}$ and  $2K_2$.
One can see that $G_1$ is $2C_k$-free   with $|G_1|=n$ and
\begin{eqnarray*}
e(G_1)&=&e(H_1)+e(2O_{k-2})+e(2K_2)+2k\\
&=&ex_{_{\mathcal{P}}}(n-2k,2C_k)+2(2(k-2)-3)+2+2k\\
&=&ex_{_{\mathcal{P}}}(n-2k,2C_k)+6k-12.
\end{eqnarray*}
Hence, $ex_{_{\mathcal{P}}}(n,2C_k)\ge e(G_1)=ex_{_{\mathcal{P}}}(n-2k,2C_k)+6k-12$, as desired.\medskip

Similarly, to prove (b),  let $H_2$ be any $C_k$-free planar graph on $n-2k$ vertices with $ex_{_{\mathcal{P}}}(n-2k,C_k)$ edges. Let $v_2\in V(H_2)$. Let $G_2$ be the planar  graph obtained from the disjoint union of $H_2\cup O_{2k-2}\cup K_2$ by joining $v_2$ to each vertex of $O_{2k-2}\cup K_2$.
Then  $G_2$ is $2C_k$-free   with $|G_2|=n$ and
\begin{eqnarray*}
e(G_2)&=&e(H_2)+e(O_{2k-2})+e(K_2)+2k\\
&=&ex_{_{\mathcal{P}}}(n-2k,C_k)+2(2k-2)-3+1+2k\\
&=&ex_{_{\mathcal{P}}}(n-2k,C_k)+6k-6.
\end{eqnarray*}
Hence, $ex_{_{\mathcal{P}}}(n,2C_k)\ge e(G_2)=ex_{_{\mathcal{P}}}(n-2k,C_k)+6k-6$, as desired.\qed
\medskip

\begin{thm}\label{CkCk+}
  Let  $n$ and  $k$  be integers with  $n\ge 2k+1$ and $k\ge4$. Then $$ex_{_{\mathcal{P}}}(n,C_k\cup C_k^+)=ex_{_{\mathcal{P}}}(n,2C_k).$$
\end{thm}

\pf   Since $2C_k$ is subgraph of $C_k\cup C_k^+$, we have $ex_{_{\mathcal{P}}}(n,C_k\cup C_k^+)\ge ex_{_{\mathcal{P}}}(n,2C_k).$ We next prove that $ex_{_{\mathcal{P}}}(n,C_k\cup C_k^+)\le ex_{_{\mathcal{P}}}(n, 2C_k)$  by induction on $n$. Let $G$ be any $C_k\cup C_k^+$-free planar graph with $|G|=n\ge2k+1$ and $e(G)=ex_{_{\mathcal{P}}}(n,C_k\cup C_k^+)$. If $G$ is $2C_k$-free, then $e(G)\le ex_{_{\mathcal{P}}}(n,2C_k)$, as desired. We may assume that $G$ contains $2C_k$ as a subgraph. Let $A$ be the set of vertices of a copy of  $2C_k$ in $G$. Then $e(G[A])\le3(2k)-6=6k-6$   and $G$ has no edges between $A$ and $V(G\less A)$ because  $k\ge4$ and $G$ is $C_k\cup C_k^+$-free.  Hence, $e(G)=e(G[A])+e(G\less A)\le 6k-6+e(G\less A)$.\medskip

We first prove that the statement holds  for $n$ satisfying   $2k+1\le n\le 4k$.  Note that $e(G\less A)\le{n-2k \choose 2}$ when $n\le2k+2$ and $e(G\less A)\le3(n-2k)-6$ when  $2k+3\le n\le 4k$. Hence,
\begin{eqnarray*}
  e(G)&\le & \begin{cases}
  6k-6 & \text{if}\,\quad  n=2k+1 \\
  6k-6+1 & \text{if}\,\quad  n=2k+2 \\
  6k-6+3(n-2k)-6 & \text{if}\,\quad 2k+3\le n\le 4k \\
\end{cases} \\
  &=&\begin{cases}
  3n-9, & \text{if}\,\quad  n=2k+1 \\
  3n-11, & \text{if}\,\quad  n=2k+2 \\
  3n-12, & \text{if}\,\quad 2k+3\le n\le 4k \\
\end{cases} \\
&\le &\left(3-\frac{1}{k-2}\right)n-5.
\end{eqnarray*}
By Lemma~\ref{22LSS}, we have
 $e(G)\le ex_{_{\mathcal{P}}}(n,2C_k)$, and so the statement holds when   $  n\le 4k$.\medskip

We may assume that $n\ge4k+1$.
We first consider the case when $G[A]$ is disconnected. Then $G[A]$ has exactly two components each of order $k$ because $G[A]$ contains a $2C_k$ as a subgraph. Then $e(G[A])\le 2(3k-6)=6k-12$. Note that   $G\less A$ is $C_k\cup C_k^+$-free and $|G\less A|=n-2k\ge 2k+1$. By the induction hypothesis, we see that   $e(G\less A)\le ex_{_{\mathcal{P}}}(n-2k,C_k\cup C_k^+) \le ex_{_{\mathcal{P}}}(n-2k,2C_k)$.
Hence,
 $e(G)=e(G[A])+e(G\less A)\le   6k-12+ex_{_{\mathcal{P}}}(n-2k,2C_k).$
By Lemma~\ref{2Ck}($a$), we see that $e(G)\le ex_{_{\mathcal{P}}}(n,2C_k)$, as desired. It remains to consider the case  $G[A]$ is connected. Then $G[A]$ contains a copy of $C_k^+$ and so $G\less A$ is $C_k$-free. Hence,  $e(G\less A)\le ex_{_{\mathcal{P}}}(n-2k,C_k)$.
It follows that
 $e(G)=e(G[A])+e(G\less A)\le 6k-6+ex_{_{\mathcal{P}}}(n-2k,C_k).$
By Lemma~\ref{2Ck}($b$), we have $e(G)\le ex_{_{\mathcal{P}}}(n,2C_k)$.\medskip

This completes the proof of Theorem~\ref{CkCk+}.\qed

\section{Planar Tur\'an numbers for  $\Theta_k^+$ and $\Theta_k$}\label{s:Theta+}
Finally,  we study the planar Tur\'an numbers for  $\Theta_k^+$ and $\Theta_k$. We prove that $ex_{_{\mathcal{P}}}(n,\Theta_k^+)=ex_{_{\mathcal{P}}}(n,\Theta_k)$ for all $n\ge k+1$ and $k\ge5$,   and $ex_{_{\mathcal{P}}}(n,\Theta_4^+)=ex_{_{\mathcal{P}}}(n,\Theta_4)$ for all $n\ge6$. We begin with a lemma.

\begin{lem}\label{Thetak}
Let $n$ and $k$ be integers with $n\ge k+2$ and $k\ge4$. Then
\begin{enumerate}[(a)]
  \item $ex_{_{\mathcal{P}}}(n,\Theta_k)\ge ex_{_{\mathcal{P}}}(n-k,\Theta_k)+3k-6$;
 \item $ex_{_{\mathcal{P}}}(n,\Theta_k)\ge ex_{_{\mathcal{P}}}(n-t,\Theta_k)+2t-1 $   for all $k\ge5$ and each $t\in\{k-1, k-2\}$; and
  \item $ex_{_{\mathcal{P}}}(n,\Theta_4)\ge ex_{_{\mathcal{P}}}(n-t,\Theta_4)+2t$  for each $t\in\{2,3\}$ and all $n\ge t+4$.
\end{enumerate}
\end{lem}

\pf   To prove $(a)$, we will construct a $\Theta_k$-free planar graph on $n$ vertices with $ex_{_{\mathcal{P}}}(n-k,\Theta_k)+3k-6$ edges. Let $H$ be any $\Theta_k$-free planar graph on  $n-k\ge 2$ vertices with $ex_{_{\mathcal{P}}}(n-k,\Theta_k)$ edges.
Let $v\in V(H)$. Let $G$ be the planar graph obtained from the disjoint copies of $H$,  $O_{k-2}$ and  $K_2$ by joining $v$ to  each vertex of $O_{k-2}$ and $ K_2$.
Then  $G$ is $\Theta_k$-free   with $|G|=n$ and
\begin{align*}
e(G)&=e(H)+e(O_{k-2})+e(K_2)+k\\
&=ex_{_{\mathcal{P}}}(n-k,\Theta_k)+(2(k-2)-3)+1+k\\
&=ex_{_{\mathcal{P}}}(n-k,\Theta_k)+3k-6.
\end{align*}
Hence, $ex_{_{\mathcal{P}}}(n,\Theta_k)\ge e(G)= ex_{_{\mathcal{P}}}(n-k,\Theta_k)+3k-6$, as desired.\medskip

 Similarly, to prove $(b)$, for  all $k\ge 5$ and each  $t\in\{k-2, k-1\}$,   let $H$ be any $\Theta_k$-free planar graph on $n-t\ge3$ vertices with $ex_{_{\mathcal{P}}}(n-t,\Theta_k)$ edges. Let $v\in V(H)$  and let $G$ be the planar graph obtained from the disjoint copies of $H$ and $ K_{1, t-1}$  by joining $v$ to  every vertex of $ K_{1, t-1}$. Then  $G$ is $\Theta_k$-free   with  $n$ vertices.   Hence, $ex_{_{\mathcal{P}}}(n,\Theta_k)\ge e(G)=e(H)+e(K_{1, t-1})+t= ex_{_{\mathcal{P}}}(n-t,\Theta_k)+2t-1$.\medskip

 It remains to prove $(c)$. For each $t\in\{2,3\}$, let $H$ be any $\Theta_4$-free plane graph on $n-t\ge 4$ vertices with $ex_{_{\mathcal{P}}}(n-t,\Theta_4)$ edges.
By Theorem~\ref{Theta4}($a$), we see that $ex_{_{\mathcal{P}}}(n-t,\Theta_4)\le {12(n-2-t)}/{5}<3(n-t)-6$. It follows that $H$ is not a plane triangulation and thus  contains a face $F$ such that $F$ is not a 3-face. Let $u_1,u_2\in V(F)$ be such that the shortest path between $u_1$ and $u_2$ has length exactly two. Let $G$ be the planar graph obtained from $H$ by adding $t$ pairwise  non-adjacent new vertices to $F$ each  joining to exactly $u_1$ and $u_2$.
Then  $G$ is $\Theta_4$-free   with $|G|=n$ and
$$e(G)=e(H)+2t=ex_{_{\mathcal{P}}}(n-t,\Theta_4)+2t.$$
Hence, $ex_{_{\mathcal{P}}}(n,\Theta_4)\ge e(G)= ex_{_{\mathcal{P}}}(n-t,\Theta_4)+2t$, as desired.\qed\medskip

\begin{thm}\label{Thetak+}
Let $n,k$ be integers with $n\ge k+1$ and $k\ge5$. Then $ex_{_\mathcal{P}}(n,\Theta_k^+)=ex_{_\mathcal{P}}(n,\Theta_k).$
\end{thm}
\pf   We   observe that  $ex_{_{\mathcal{P}}}(n,\Theta_k^+)\ge ex_{_{\mathcal{P}}}(n,\Theta_k)$,  because for each  $G\in \Theta_k$, there exists  $G'\in \Theta_k^+$ such that $G$ is a subgraph of $G'$. We next prove that $ ex_{_{\mathcal{P}}}(n,\Theta_k^+)\le ex_{_{\mathcal{P}}}(n,\Theta_k)$. Suppose not. Choose a minimum $n$ such that $n\ge k+1$ and  $ ex_{_{\mathcal{P}}}(n,\Theta_k^+)> ex_{_{\mathcal{P}}}(n,\Theta_k)$.   Let $G$ be a $\Theta_k^+$-free planar graph with  $|G|=n$  and $e(G)=ex_{_{\mathcal{P}}}(n,\Theta_k^+)>ex_{_{\mathcal{P}}}(n,\Theta_k)$.    If $G$ is $\Theta_k$-free, then $e(G)\le ex_{_{\mathcal{P}}}(n,\Theta_k)$,  a contradiction. Thus  $G$ contains some $H^*\in\Theta_k$ as a subgraph. Let $u_1,u_2$ be the vertices of degree three in $H^*$, and let $H:=G[V(H^*)]$ and $X:=V(G)\less V(H)$.
Since $G$ is $\Theta_k^+$-free, we see that $N_G(x)\cap  V(H)\subseteq \{u_1,u_2\}$ for all $x\in X$. Let $U\subseteq \{u_1, u_2\}$ such that each $u\in U$ is adjacent to some vertex in $X$. Then $0\le |U|\le 2$. We may assume that $u_1\in U$ when $U\ne\emptyset$. Note that  when $U=\{u_1, u_2\}$, then $u_1u_2$  is an edge of both  $H$ and $G[X\cup U]$. This implies that $e(G)=e(H)+e(G[X\cup U])-\max\{0, |U|-1\}$. We claim that
 \begin{align*}
  e(H)\le & \begin{cases}
  3k-6 & \text{if}\,\quad  |U|=0 \\
  2k-3 & \text{if}\,\quad  |U|=1 \\
  k+1 & \text{if}\,\quad |U|=2.\\
\end{cases}
\end{align*}
When $|U|=0$, then   $H$ is a component of $G$ and so $e(H)\le 3k-6$. Next when $|U|=1$, let $w_1, w_2$ be the other two neighbors of $u_1$ in $H^*$. Then $d_H(u_2)=3$, $d_H(v)=2$ for each $v\in \{w_1, w_2\}$,  and  $d_{H\less u_1}(v)= 2$ for all $v\in  V(H)\less \{u_1, u_2, w_1, w_2\}$, 
otherwise $G$ is not $\Theta_k^+$-free. Thus $H$ is an outer-planar graph and so $e(H)\le 2k-3$. Finally, when $U=\{u_1, u_2\}$, then  $d_H(u_1)=d_H(u_2)=3$ and $d_H(v)=2$ for all $v\in V(H)\less\{u_1, u_2\}$, otherwise $G$ is not $\Theta_k^+$-free. Hence $e(H)=k+1$. This proves the claim and so
\begin{align*}\tag{$*$}
e(G)&\le    e(H)+e(G[X\cup U])-\max\{0, |U|-1\}\le \begin{cases}
  3k-6 +e(G[X\cup U])& \text{if}\,\quad  |U|=0 \\
  2k-3+e(G[X\cup U]) & \text{if}\,\quad  |U|=1 \\
  k+e(G[X\cup U])& \text{if}\,\quad |U|=2.\\
  \end{cases} \medskip
\end{align*}

Suppose $n\ge 2k-|U|+1$. Then $|X\cup U|=n-k+|U|\ge k+1$ and $G[X\cup U]$ is $\Theta_k^+$-free. Thus  $e(G[X\cup U])\le ex_{_\mathcal{P}}(n-k+|U|,\Theta_k^+)$. By the minimality of $n$,  $ex_{_\mathcal{P}}(n-k+|U|,\Theta_k^+)\le ex_{_\mathcal{P}}(n-k+|U|,\Theta_k)$. By ($*$), we have
\begin{align*}
e(G) \le \begin{cases}
  3k-6 +ex_{_\mathcal{P}}(n-k,\Theta_k)& \text{if}\,\quad  |U|=0 \\
  2k-3+ex_{_\mathcal{P}}(n-k+1,\Theta_k) & \text{if}\,\quad  |U|=1 \\
  k+ex_{_\mathcal{P}}(n-k+2,\Theta_k)& \text{if}\,\quad |U|=2.\\
  \end{cases} \medskip
\end{align*}
By Lemma~\ref{Thetak}($a,b$),  we see that $e(G)\le ex_{_\mathcal{P}}(n,\Theta_k)$, a contradiction. This proves that $k+1\le n\le 2k-|U|$.  Let $G^*$ be the graph defined as in the proof of Theorem~\ref{Ck+}.   Note that  $G^*$ is   $\Theta_k$-free. It suffices to show that    $e(G)\le e(G^*)$ so we obtain a desired contradiction. \medskip

We first consider the case $|U|=0$. Then $k+1\le n\le 2k$ and $1\le |X|=n-k\le k$.  Thus  $e(G[X])\le {|X|\choose 2}$ when $k+1\le n\le k+2$, and $e(G[X])\le 3|X|-6$ when $k+3\le n\le 2k$. By ($*$),
\begin{eqnarray*}
  e(G)\le  3k-6+e(G[X])&\le & \begin{cases}
  3k-6 & \text{if}\,\quad  n=k+1 \\
  3k-6+1 & \text{if}\,\quad  n=k+2 \\
  3k-6+3(n-k)-6 & \text{if}\,\quad k+3\le n\le 2k \\
\end{cases} \\
  &=&\begin{cases}
   3n-9 & \text{if}\,\quad  n=k+1 \\
  3n-11 & \text{if}\,\quad  n=k+2 \\
  3n-12 & \text{if}\,\quad k+3\le n\le 2k \\
\end{cases}\\
&\le & e(G^*).
\end{eqnarray*}
 We next  consider the case $|U|=1$. Then $k+1\le n\le 2k-1$ and $2\le |X\cup U|=n-k+1\le k$.  Note that when $k=5$ and $n=2k-1$,  we have $|X\cup U|=5$ and $G[X\cup U]$ is not the plane triangulation on $5$ vertices, else   $G[X\cup \{u_1, u_2\}]$  contains a copy of $\Theta_5^+$.
Thus  $e(G[X\cup U])\le 3|X\cup U|-7$ when $n=2k-1$ and $k=5$;  $e(G[X\cup U])\le 1$ when $n=k+1$;  and $e(G[X\cup U])\le 3|X\cup U|-6$ when $k+2\le n\le 2k-2$  or $n=2k-1$ and $k\ge6$. By ($*$),
\begin{eqnarray*}
  e(G)\le 2k-3+e(G[X\cup U])
 & \le&  \begin{cases}
  2k-3+1 & \text{if}\,\quad  n=k+1 \\
  2k-3+3(n-k+1)-6 & \text{if}\,\quad k+2\le n\le 2k-2\\
  2k-3+3(n-k+1)-6 & \text{if}\,\quad n=2k-1\, \text{and}\,\, k\ge6\\
  2k-3+3(n-k+1)-7 & \text{if}\,\quad n=2k-1\, \text{and}\,\, k=5\\
\end{cases} \\
  &=&\begin{cases}
   3n-k-5 & \text{if}\,\quad  n=k+1 \\
  3n-k-6 & \text{if}\,\quad k+2\le n\le 2k-2 \\
    3n-k-6 & \text{if}\,\quad n=2k-1\, \text{and}\, k\ge6\\
  3n-k-7 & \text{if}\,\quad n=2k-1\, \text{and}\, k=5\\
\end{cases}\\
&\le &e(G^*).
\end{eqnarray*}
  Finally, suppose  $|U|=2$. 
  Then $k+1\le n\le 2k-2$ and  $3\le |X\cup U|=n-k+2\le k$.  Note that when $k=5$ and $n=2k-2$,  we have $|X\cup U|=5$ and $G[X\cup U]$ is not the plane triangulation on $5$ vertices, else   $G$  contains a copy of $\Theta_5^+$.  Thus   $e(G[X\cup U])\le 3|X\cup U|-7$ when $n=2k-2$ and   $  k=5$; and $e(G[X\cup U])\le 3|X\cup U|-6$ when $k+1\le n\le 2k-3$   or $n=2k-2$ and $k\ge6$.  By ($*$),
\begin{eqnarray*}
  e(G)\le  k+e(G[X\cup U])
 & \le&  \begin{cases}
  k+3(n-k+2)-6 & \text{if}\,\quad k+1\le n\le 2k-3\\
  k+3(n-k+2)-6 & \text{if}\,\quad n=2k-2\, \text{and}\, k\ge6\\
  k+3(n-k+2)-7 & \text{if}\,\quad n=2k-2\, \text{and}\, k=5\\
\end{cases} \\
  &=&\begin{cases}
  3n-2k & \text{if}\,\quad k+1\le n\le 2k-3 \\
    3n-2k & \text{if}\,\quad n=2k-2\, \text{and}\, k\ge6\\
  3n-2k-1 & \text{if}\,\quad n=2k-2\, \text{and}\, k=5\\
\end{cases}\\
&\le & e(G^*).
\end{eqnarray*}
This completes the proof of Theorem~\ref{Thetak+}.\qed\bigskip

\begin{thm}\label{Theta4+}
  Let $n\ge5$ be integer. Then
  \[ex_{_\mathcal{P}}(n,\Theta_4^+)=
  \begin{cases}
    ex_{_\mathcal{P}}(n,\Theta_4) & \mbox{if } n\ge 6 \\
    ex_{_\mathcal{P}}(n,\Theta_4)+1=7 & \mbox{if } n=5.
  \end{cases}  \]
\end{thm}

\pf  Throughout the proof, by abusing notation, we use $\Theta_4$ and $\Theta_4^+$ to denote the unique graph in $\Theta_4$ and $\Theta_4^+$, respectively. We first observe that $ex_{_\mathcal{P}}(5,\Theta_4)=6$ because  $K_{2,3}$ is $\Theta_4$-free; and every planar graph with five vertices and seven edges is  obtained from $K_5$ with three edges removed,  and it is simple to check that   every such graph obtained from $K_5$     contains a copy of $\Theta_4$. Similarly,   $ex_{_\mathcal{P}}(5,\Theta_4^+)=7$ because  $K_2+\overline{K}_3$ is $\Theta_4^+$-free; and every planar graph with five vertices and eight edges is  obtained from $K_5$ with two edges removed, and   every such graph  obtained from $K_5$  contains a copy of   $\Theta_4^+$.  This proves that $ex_{_\mathcal{P}}(5,\Theta_4^+)=ex_{_\mathcal{P}}(5,\Theta_4)+1=7$.  We next prove that $ex_{_\mathcal{P}}(n,\Theta_4^+)=ex_{_\mathcal{P}}(n,\Theta_4)$ for all $n\ge 6$.
\medskip

Since $\Theta_4$ is a subgraph of $\Theta_4^+$, we see that  $ex_{_{\mathcal{P}}}(n,\Theta_4^+)\ge ex_{_{\mathcal{P}}}(n,\Theta_4)$.  We next   prove that $ ex_{_{\mathcal{P}}}(n,\Theta_4^+)\le ex_{_{\mathcal{P}}}(n,\Theta_4)$ for all $n\ge 6$ by induction on $n$. Let $G$ be a $\Theta_4^+$-free planar graph on $n\ge 6$ vertices with $ex_{_{\mathcal{P}}}(n,\Theta_4^+)$ edges. If $G$ is $\Theta_4$-free, then $e(G)\le ex_{_{\mathcal{P}}}(n,\Theta_4)$, as desired. We may assume that $G$ contains $\Theta_4$ as a subgraph.  Let $u_1,u_2, u_3, u_4$ be the vertices of  a copy of  $H^*:=\Theta_4$ in $G$ such that $d_{H^*}(u_1)=d_{H^*}(u_2)=3$, and let $H:=G[V(H^*)]$ and $X:=V(G)\less V(H)$.
Since $G$ is $\Theta_4^+$-free, we see that $N_G(x)\cap  V(H)\subseteq \{u_1,u_2\}$ for all $x\in X$. Let $U\subseteq \{u_1, u_2\}$ be such that each $u\in U$ is adjacent to some vertex in $X$. Then $0\le |U|\le 2$. Note that  $d_H(u_3)=d_H(u_4)=2$ when $U\ne \emptyset$, and $u_1u_2$  is an edge of both  $H$ and $G[X\cup U]$ when $U=\{u_1, u_2\}$. Thus $e(H)=5+\max\{0, 1-|U|\}$ and so \[e(G)=e(H)+e(G[X\cup U])-\max\{0,|U|-1\}=6-|U|+e(G[X\cup U]).\tag{$\dagger$}\]

\begin{figure}[htbp]
\centering
\includegraphics*[width=4cm]{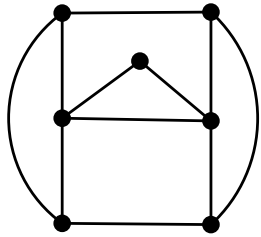}
\caption{Graph $J$ when $n=7$ and $|U|=2$.}\label{J}
\end{figure}

  For $6\le n\le 9-|U|$, we have $2\le |X\cup U|=n-4+|U|\le 5$. We first consider the case  $|X\cup U|=5$. Then $n=9-|U|\ge7$ and $e(G[X\cup U])\le ex_{_\mathcal{P}}(5,\Theta_4^+)=7$ and so $e(G)\le 6-|U|+7=13-|U|$ by ($\dagger$).  Let $J:=K_{2, 7-|U|}$ when $|U|\le 1$ and let $J$ be the graph  given in  Figure~\ref{J} when $|U|=2$. Then $J$ is a $\Theta_4$-free planar graph with $n=9-|U|$ vertices, and so
  $e(G)\le e(J)\le ex_{_\mathcal{P}}(9-|U|,\Theta_4)$.  We next consider the case $|X\cup U|=4$. Then $n=8-|U|\ge6$; in addition, if $U\ne \emptyset$, then $G[X\cup U]\ne K_4$, else $G$ contains a copy of $\Theta_4^+$.  This implies that $e(G[X\cup U])\le 6-\min\{1, |U|\}$ and so $e(G)\le 6-|U|+6-\min\{1, |U|\}=12-|U|-\min\{1, |U|\}$ by ($\dagger$).  Let $J:=K_{2, 6-|U|}$ when $|U|\le 1$ and $J:=\overline{C_6}$  when $|U|=2$. Then $J$ is a $\Theta_4$-free planar graph with $n=8-|U|$ vertices.
It follows that $e(G)\le e(J)\le ex_{_\mathcal{P}}(8-|U|,\Theta_4)$. Hence $2\le |X\cup U|\le3$. In this case,  $n=4+|X\cup U|-|U|\le 7-|U|$  and  $|U|\le 1$ because $n\ge6$.
In addition,  $e(G[X\cup U])\le {|X\cup U|\choose 2}$. By ($\dagger$),
\[e(G) \le  6-|U|+ {|X\cup U|\choose 2}\le 2(2+|X\cup U|-|U|)=e(K_{2,  n-2})\le ex_{_\mathcal{P}}(n,\Theta_4),\]
because $K_{2,  n-2}$ is $\Theta_4$-free. \medskip

We may assume that  $n\ge 10-|U|$. Then $|X\cup U|=n-4+|U|\ge 6$ and $G[X\cup U]$ is $\Theta_4^+$-free. Thus  $e(G[X\cup U])\le ex_{_\mathcal{P}}(n-4+|U|,\Theta_4^+)$. By the  induction hypothesis,  $ex_{_\mathcal{P}}(n-4+|U|,\Theta_4^+)\le ex_{_\mathcal{P}}(n-4+|U|,\Theta_4)$. By ($\dagger$),
\begin{align*}
e(G)\le   6-|U|+ex_{_\mathcal{P}}(n-4+|U|,\Theta_4^+)
 \le \begin{cases}
  6 +ex_{_\mathcal{P}}(n-k,\Theta_k)& \text{if}\,\quad  |U|=0 \\
5+ex_{_\mathcal{P}}(n-k+1,\Theta_k) & \text{if}\,\quad  |U|=1 \\
4 +ex_{_\mathcal{P}}(n-k+2,\Theta_k)& \text{if}\,\quad |U|=2.\\
  \end{cases} \medskip
\end{align*}
By Lemma~\ref{Thetak}($a,c$),  we see that $e(G)\le ex_{_\mathcal{P}}(n,\Theta_4)$, as desired.
 \qed\bigskip

Corollary~\ref{c:Thetak+}  follows directly from Theorem~\ref{Thetak+},  Theorem~\ref{Theta4+}, Theorem~\ref{Theta4}($a,b$), Theorem~\ref{C6Theta6}($b$) and Theorem~\ref{LB-Ck}.

\begin{cor}\label{c:Thetak+}
Let $n$ be a positive integer. Then

 \begin{enumerate}[(a)]

\item $ex_{_\mathcal{P}}(n, \Theta_4^+)\leq     {12(n-2)}/5$ for all $n\ge5$, with  equality when $n\equiv 12 (\rm{mod}\, 20)$.

\item  $ex_{_\mathcal{P}}(n, \Theta_5^+)\leq {5(n-2)}/2$ for all $n\ge6$, with  equality when $n\equiv 50 (\rm{mod}\, 120)$.

\item  $ex_{_\mathcal{P}}(n, \Theta_6^+)\leq(18n-48)/7$ for all $n\geq 14$. Equality  holds for infinitely many $n$.

\item $ex_{_{\mathcal{P}}}(n,\Theta_k^+)=3n-6$ for all $k\ge 11$ and $k+1\le n\le  k-5+\lfloor(k-1)/2\rfloor$, and   \[ex_{_{\mathcal{P}}}(n,\Theta_k^+)\ge \left(3-\frac{3-\frac{2}{k-1}}{k-6+\left\lfloor\frac{k-1}2\right\rfloor}\right)n+
  \frac{12+3r-\frac{8+2r}{k-1}}{k-6+\left\lfloor\frac{k-1}2\right\rfloor}+\frac{4}{k-1}-\min\{r+10,11\}\] for all $k\ge11$ and  $n\ge  k-4+\lfloor{(k-1)}/{2}\rfloor$, where $r$ is the remainder of  $n-4$ when divided by $k-6+\lfloor(k-1)/2\rfloor$.

\end{enumerate}
\end{cor}


\frenchspacing
\begin{thebibliography}{14}
\bibitem{CLLS21}
D. Cranston, B. Lidick\'y, X. Liu,  A. Shantanam, Planar Tur\'an numbers of cycles: a counterexample,  Electron. J. Combin. 29(3) (2022), \#P3.31.

\bibitem{Dowden2016}
C. Dowden, Extremal $C_4$-free/$C_5$-free planar graphs,  J. Graph Theory  83 (2016),  213--230.


\bibitem{ErdosStone} P.  Erd\H{o}s,  A. Stone, On the structure of linear graphs,  Bull. Amer. Math. Soc. 52 (1946), 1087--1091.

\bibitem{LZW}
L. Fang, B. Wang,  M. Zhai, Planar Tur\'an number of intersecting triangles, Discrete Math.  345(5) (2022), 112794.


\bibitem{C6}
D. Ghosh, E. Gy\H{o}ri, R. Martin, A. Paulos,  C. Xiao, Planar Tur\'an number of the 6-cycle, SIAM J. Discrete Math. 36(3) (2022), 2028--2050.

\bibitem{Theta6}
D. Ghosh, E. Gy\H{o}ri,   A. Paulos, C. Xiao,  O. Zamora, Planar Tur\'an number of the $\Theta_6$-cycle,  arXiv:2006.00994v1.

\bibitem{doublestar}
D. Ghosh, E. Gy\H{o}ri,   A. Paulos,  C. Xiao, Planar Tur\'an number of double stars,  arXiv:2110.10515v2.

\bibitem{15HJSS}
M. Hor\v{n}\'ak, S. Jendrol$'$, I. Schiermeyer, R. Sot\'ak, Rainbow numbers for cycles in plane triangulations, J. Graph Theory 78 (2015), 248--257.

\bibitem{19LSSan}
Y. Lan, Y. Shi and Z.-X. Song, Planar anti-Ramsey numbers of paths and cycles, Discrete Math. 342 (2019), 3216--3224.

\bibitem{LSS}
Y. Lan, Y. Shi and Z.-X. Song, Extremal Theta-free planar graphs, Discrete Math. 342 (2019), 111610.


\bibitem{19LSS}
Y. Lan, Y. Shi and Z.-X. Song, Extremal $H$-free planar graphs, Electron. J. Combin. 26 (2019), \#P2.11.


\bibitem{LSS21}
Y. Lan, Y. Shi, Z.-X. Song,  Planar Tur\'an number and planar anti-Ramsey number of graphs, Oper. Res. Trans.  25 (2021), 200--216.

\bibitem{22LSS}
Y. Lan, Y. Shi and Z.-X. Song,
Planar Tur\'an numbers of cubic graphs and disjoint union of cycles, arXiv:2202.09216v2.



\bibitem{Turan1941}
P.  Tur\'an, On an extremal problem in graph theory,  Mat.  Fiz. Lapok 48 (1941), 436--452.

 \end{thebibliography}
\end{document}